\documentclass[12pt]{article}
\usepackage{amssymb, amsmath, url, graphicx}

\def\3{\subset }
\def\4{\subseteq }
\def\<{\left<}
\def\>{\right>}

\def\bit{\begin{itemize}}
\def\eit{\end{itemize}}
\def\3{\subset }
\def\4{\subseteq }

\def\0{\leqno}

\def\barr{\begin{array}}
\def\earr{\end{array}}
\def\dd{\displaystyle}

\def\Z{{\rlap{$\kern2pt{\rm Z}$}{\rm Z}\,}}

\title{\bf Factorization numbers of finite rank $3$ abelian $p$-groups}
\author{Marius T\u arn\u auceanu}
\date{March 22, 2016}

\begin{document}

\maketitle

\begin{abstract}
In this short note we give a formula for the factorization number $F_2(G)$ of
a finite rank $3$ abelian $p$-group $G$. This extends a result in our previous
work \cite{9}.
\end{abstract}

\noindent{\bf MSC (2010):} Primary 20D40; Secondary 20D60.

\noindent{\bf Key words:} factorization number, subgroup
commutativity degree, M\"{o}bius function, finite abelian $p$-group.

\section{Introduction}

Let $G$ be a finite group, $L(G)$ be the subgroup lattice of $G$
and $H$, $K$ be two subgroups of $G$. If $G=HK$, then $G$ is said
to be \textit{factorized} by $H$ and $K$ and the expression $G=HK$
is said to be a \textit{factorization} of $G$. Denote by $F_2(G)$
the \textit{factorization number} of $G$, that is the number of all
factorizations of $G$.

The starting point for our discussion is given by the papers
\cite{6,1}, where $F_2(G)$ has been computed for certain classes of
finite groups. Then, in \cite{9}, we have obtained explicit formulas for the
factorization numbers of an elementary abelian $p$-group and of a rank 2
abelian $p$-group. These are based on the connection between $F_2(G)$ and
the subgroup commutativity degree $sd(G)$ of $G$ (see \cite{7,8}), namely
$$sd(G)=\dd\frac{1}{|L(G)|^2}\sum_{H\leq G} F_2(H).$$Obviously, by
applying the well-known M\"{o}bius inversion formula to the above
equality, one obtains
$$F_2(G)=\dd\sum_{H\leq G} sd(H)|L(H)|^2\mu(H,G).\0(1)$$In
particular, if $G$ is abelian, then we have $sd(H)=1$ for all
$H\in L(G)$, and consequently
$$F_2(G)=\dd\sum_{H\leq G} |L(H)|^2\mu(H,G)=\dd\sum_{H\leq G} |L(G/H)|^2\mu(H).\0(2)$$This
formula will be used in the following to calculate the factorization number of a rank $3$ abelian
$p$-group.
\bigskip

First of all, we recall a theorem due to P. Hall \cite{2} (see
also \cite{4}), that permits us to compute explicitly the
M\"{o}bius function of a finite $p$-group.

\bigskip\noindent{\bf Theorem 1.} {\it Let $G$ be a finite $p$-group of order $p^n$.
Then $\mu(G)=0$ unless $G$ is elementary abelian, in which case we
have $\mu(G)=(-1)^n p^{\,\binom{n}{2}}.$}
\bigskip

We also need to know the total number of subgroups of a finite rank $3$ abelian
$p$-group. It has been determined by using different methods in \cite{3,5,10}.

\bigskip\noindent{\bf Theorem 2.} {\it The total number of subgroups of $\mathbb{Z}_{p^{\lambda_1}}{\times}\mathbb{Z}_{p^{\lambda_2}}{\times}\mathbb{Z}_{p^{\lambda_3}}$,
where $\lambda_1\geq\lambda_2\geq \lambda_3\geq 1$, is
\begin{equation*}
f(\lambda_1,\lambda_2,\lambda_3)=\frac{A}{(p^2-1)^2(p-1)}\,,\0(3)
\end{equation*}
where
\begin{align*}
A=&(\lambda_3{+}1)(\lambda_1{-}\lambda_2{+}1)p^{\lambda_2{+}\lambda_3{+}5}{+}2(\lambda_3{+}1)p^{\lambda_2{+}\lambda_3{+}4}{-}2(\lambda_3{+}1)(\lambda_1{-}\lambda_2)p^{\lambda_2{+}\lambda_3{+}3}\\
&{-}2(\lambda_3{+}1)p^{\lambda_2{+}\lambda_3{+}2}{+}(\lambda_3{+}1)(\lambda_1{-}\lambda_2{-}1)p^{\lambda_2{+}\lambda_3{+}1}{-}(\lambda_1{+}\lambda_2{-}\lambda_3{+}3)p^{2\lambda_3{+}4}\\
&{-}2 p^{2\lambda_3{+}3}{+}(\lambda_1{+}\lambda_2{-}\lambda_3{-}1) p^{2\lambda_3{+}2}{+}(\lambda_1{+}\lambda_2{+}\lambda_3{+}5)p^2{+}2p{-}(\lambda_1{+}\lambda_2{+}\lambda_3{+}1).
\end{align*}}

We are now able to give the main result of our note.

\bigskip\noindent{\bf Theorem 3.} {\it The following equality holds
\begin{align*}
F_2(\mathbb{Z}_{p^{\lambda_1}}{\times}\mathbb{Z}_{p^{\lambda_2}}{\times}\mathbb{Z}_{p^{\lambda_3}}){=}&{-}p^3f^2(\lambda_1{-}1,\lambda_2{-}1,\lambda_3{-}1){+}p(f^2(\lambda_1{-}1,\lambda_2{-}1,\lambda_3)\\
&{+}pf^2(\lambda_1{-}1,\lambda_2,\lambda_3{-}1){+}p^2f^2(\lambda_1,\lambda_2{-}1,\lambda_3{-}1))\\
&{-}(f^2(\lambda_1{-}1,\lambda_2,\lambda_3){+}pf^2(\lambda_1,\lambda_2{-}1,\lambda_3)\\
&{+}p^2f^2(\lambda_1,\lambda_2,\lambda_3{-}1)){+}f^2(\lambda_1,\lambda_2,\lambda_3),
\end{align*}
where the quantities $f(\lambda_1,\lambda_2,\lambda_3)$ are given by $(3)$.}
\bigskip

We remark that for $\lambda_3=0$ the above equality leads to the formula in Theorem 3 of \cite{9}. It simplifies in the particular case $\lambda_1=\lambda_2=\lambda_3=\lambda$.

\bigskip\noindent{\bf Corollary 4.} {\it We have
\begin{align*}
F_2(\mathbb{Z}_{p^{\lambda}}{\times}\mathbb{Z}_{p^{\lambda}}{\times}\mathbb{Z}_{p^{\lambda}}){=}&{-}p^3f^2(\lambda{-}1,\lambda{-}1,\lambda{-}1){+}p(1{+}p{+}p^2)f^2(\lambda,\lambda{-}1,\lambda{-}1)\\
&{-}(1{+}p{+}p^2)f^2(\lambda,\lambda,\lambda{-}1){+}f^2(\lambda,\lambda,\lambda).
\end{align*}}

However, even in this case, an explicit formula for the factorization number is too difficult to be written, but we can do it for small values of $\lambda_1,\lambda_2,\lambda_3$.

\bigskip\noindent{\bf Examples.}
\begin{itemize}
\item[a)] $F_2(\mathbb{Z}_{p^3}{\times}\mathbb{Z}_{p^2}{\times}\mathbb{Z}_p)=9p^6+15p^5+21p^4+16p^3+20p^2+11p+13$.
\item[b)] $F_2(\mathbb{Z}_{p^2}{\times}\mathbb{Z}_{p^2}{\times}\mathbb{Z}_{p^2})=5p^8+7p^7+16p^6+15p^5+21p^4+16p^3+20p^2+11p+13$.
\end{itemize}

\section{Proof of Theorem 3}

It is well-known that $G=\mathbb{Z}_{p^{\lambda_1}}{\times}\mathbb{Z}_{p^{\lambda_2}}{\times}\mathbb{Z}_{p^{\lambda_3}}$
has a unique elementary abelian subgroup of order $p^3$, say $M$,
and that
$$G/M\cong\Phi(G)\cong\mathbb{Z}_{p^{\lambda_1-1}}{\times}\mathbb{Z}_{p^{\lambda_2-1}}{\times}\mathbb{Z}_{p^{\lambda_3-1}},$$where $\Phi(G)$ 
is the Frattini subgroup of $G$. Moreover,
all elementary abelian subgroups of $G$ are contained in
$M$. Denote by $M_i$, $i=1,2,...,p^2+p+1$, the subgroups of order $p$ and by $M'_i$, $i=1,2,...,p^2+p+1$, the subgroups of order $p^2$ in
$M$. Then every quotient $G/M_i$ is isomorphic to a maximal
subgroup of $G$, which are: $p^2$ of type $(\lambda_1,\lambda_2,\lambda_3-1)$, $p$ of type $(\lambda_1,\lambda_2-1,\lambda_3)$, and $1$ of type $(\lambda_1-1,\lambda_2,\lambda_3)$. Similarly, every quotient $G/M'_i$ is isomorphic to a subgroup of index $p^2$ of $G$ that contains $\Phi(G)$, and these are: $p^2$ of type $(\lambda_1,\lambda_2-1,\lambda_3-1)$, $p$ of type $(\lambda_1-1,\lambda_2,\lambda_3-1)$, and $1$ of type $(\lambda_1-1,\lambda_2-1,\lambda_3)$. In this way the equality (2) becomes
\begin{align*}
F_2(G)=&|L(G/M)|^2\mu(M)+\dd\sum_{i=1}^{p^2+p+1}|L(G/M_i)|^2\mu(M_i)\\
&+\dd\sum_{i=1}^{p^2+p+1}|L(G/M'_i)|^2\mu(M'_i)+|L(G)|^2\mu(1),
\end{align*}
in view of Theorem 1. Since we have $\mu(M)=\mu(\mathbb{Z}_p^3)=-p^3$,
$\mu(M_i)=\mu(\mathbb{Z}_p)=-1$, $\mu(M'_i)=\mu(\mathbb{Z}_p^2)=p$, $\forall\, i=1,2,...,p^2+p+1$, and
$\mu(1)=1$, one obtains
\begin{align*}
F_2(G){=}&{-}p^3|L(\mathbb{Z}_{p^{\lambda_1-1}}{\times}\mathbb{Z}_{p^{\lambda_2-1}}{\times}\mathbb{Z}_{p^{\lambda_3-1}})|^2{+}p(|L(\mathbb{Z}_{p^{\lambda_1-1}}{\times}\mathbb{Z}_{p^{\lambda_2-1}}{\times}\mathbb{Z}_{p^{\lambda_3}})|^2\\
&{+}p|L(\mathbb{Z}_{p^{\lambda_1-1}}{\times}\mathbb{Z}_{p^{\lambda_2}}{\times}\mathbb{Z}_{p^{\lambda_3-1}})|^2{+}p^2|L(\mathbb{Z}_{p^{\lambda_1}}{\times}\mathbb{Z}_{p^{\lambda_2-1}}{\times}\mathbb{Z}_{p^{\lambda_3-1}})|^2)\\
&{-}(|L(\mathbb{Z}_{p^{\lambda_1-1}}{\times}\mathbb{Z}_{p^{\lambda_2}}{\times}\mathbb{Z}_{p^{\lambda_3}})|^2{+}p|L(\mathbb{Z}_{p^{\lambda_1}}{\times}\mathbb{Z}_{p^{\lambda_2-1}}{\times}\mathbb{Z}_{p^{\lambda_3}})|^2\\
&{+}p^2|L(\mathbb{Z}_{p^{\lambda_1}}{\times}\mathbb{Z}_{p^{\lambda_2}}{\times}\mathbb{Z}_{p^{\lambda_3-1}})|^2){+}|L(\mathbb{Z}_{p^{\lambda_1}}{\times}\mathbb{Z}_{p^{\lambda_2}}{\times}\mathbb{Z}_{p^{\lambda_3}})|^2.
\end{align*}
Under the notation of Theorem 2 this leads to the desired formula.
\hfill\rule{1,5mm}{1,5mm}

\vspace*{5ex}\small

\hfill
\begin{minipage}[t]{5cm}
Marius T\u arn\u auceanu \\
Faculty of  Mathematics \\
``Al.I. Cuza'' University \\
Ia\c si, Romania \\
e-mail: {\tt tarnauc@uaic.ro}
\end{minipage}

\end{document}